\documentstyle[11pt]{article}

\def \qed {\hspace*{\fill}\frame{\rule[0pt]{0pt}{8pt}\rule[0pt]{8pt}{0pt}}\par}
\def\F{{\rm I\!F}}
\def\P{{\rm I\!P}}
\def\C{{\mathchoice {\setbox0=\hbox{$\displaystyle\rm C$}\hbox{\hbox
to0pt{\kern0.4\wd0\vrule height0.9\ht0\hss}\box0}}
{\setbox0=\hbox{$\textstyle\rm C$}\hbox{\hbox
to0pt{\kern0.4\wd0\vrule height0.9\ht0\hss}\box0}}
{\setbox0=\hbox{$\scriptstyle\rm C$}\hbox{\hbox
to0pt{\kern0.4\wd0\vrule height0.9\ht0\hss}\box0}}
{\setbox0=\hbox{$\scriptscriptstyle\rm C$}\hbox{\hbox
to0pt{\kern0.4\wd0\vrule height0.9\ht0\hss}\box0}}}}

\def\Q{{\mathchoice {\setbox0=\hbox{$\displaystyle\rm Q$}\hbox{\raise
0.15\ht0\hbox to0pt{\kern0.4\wd0\vrule height0.8\ht0\hss}\box0}}
{\setbox0=\hbox{$\textstyle\rm Q$}\hbox{\raise
0.15\ht0\hbox to0pt{\kern0.4\wd0\vrule height0.8\ht0\hss}\box0}}
{\setbox0=\hbox{$\scriptstyle\rm Q$}\hbox{\raise
0.15\ht0\hbox to0pt{\kern0.4\wd0\vrule height0.7\ht0\hss}\box0}}
{\setbox0=\hbox{$\scriptscriptstyle\rm Q$}\hbox{\raise
0.15\ht0\hbox to0pt{\kern0.4\wd0\vrule height0.7\ht0\hss}\box0}}}}

\def\Z{{\mathchoice {\hbox{$\sf\textstyle Z\kern-0.4em Z$}}
{\hbox{$\sf\textstyle Z\kern-0.4em Z$}}
{\hbox{$\sf\scriptstyle Z\kern-0.3em Z$}}
{\hbox{$\sf\scriptscriptstyle Z\kern-0.2em Z$}}}}

\newcommand{\D}{\displaystyle} 
\newcommand{\vv}{\vspace{0.3cm}} 
\newcommand{\vV}{\vspace{0.5cm}}

\begin{document}
\title{K3 Surfaces with Nine Cusps}
\author{W. Barth} 
\maketitle 

\begin{center}
Mathematisches Institut der Universit\"{a}t

Bismarckstr. 1 1/2, D 91054 Erlangen
\end{center}

\vspace{0.5cm}
\begin{abstract}
By a $K3$-surface with nine cusps I mean a surface with nine
isolated double points $A_2$, but otherwise smooth, such that
its minimal desingularisation is a $K3$-surface. It is shown, that
such a surface admits a cyclic triple cover branched precisely
over the cusps. This parallels the theorem of Nikulin,
that a $K3$-surface with 16 nodes is a Kummer quotient
of a complex torus.
\end{abstract}

\vspace{0.5cm}
\tableofcontents

\section{Introduction}
If $E_1,...,E_{16}$ are 16 disjoint, smooth curves on a $K3$-surface
$X$ then the divisor $\sum_1^{16} E_i$ is divisible by $2$ in
$Pic(X)$. This was observed by V.V. Nikulin [N]. Equivalently: If 
$\bar{X}$ is the surface obtained from $X$ by blowing
down the 16 rational curves to nodes $e_i \in \bar{X}$,
there is a double cover $A \to \bar{X}$, with $A$ a complex torus, 
branched exactly over the 16 nodes $e_i$. The
surface $\bar{X}$ is the Kummer surface of the 
complex torus $A$.

The aim of this note is to prove an analog of Nikulin's theorem
in the case of nine cusps (double points $A_2$) instead of 16 nodes 
(double points $A_1$): 

{\bf Theorem.} {\em Let $E_i,E_i', i=1,...,9$ be 18
smooth rational curves on a $K3$-surface $X$ with
$$E_i.E_i'=1, \quad E_i.E_j=E_i.E_j'=E_i'.E_j'=0 \mbox{ for }
i \not=j,$$
then there are integers $a_i,a_i'=1,2, a_i \not=a_i',$ such that
the divisor $\sum_1^9 (a_iE_i+a_i'E_i')$ is divisible by $3$
in $Pic(X)$. Equivalently: If $\bar{X}$ is the surface
obtained from $X$ by blowing down the nine pairs of
rational curves to cusps $e_i \in \bar{X}$, then there
is a cyclic cover $A \to \bar{X}$ of order three, with
$A$ a complex torus, branched exactly over the
nine cusps $e_i$.}

The proof I give here essentially parallels Nikulin's
proof in [N].

In the case of Nikulin's theorem of course each complex
torus $A$ of dimension two appears (the covering involution
is the map $a \mapsto -a$). But complex tori of dimension
two admitting an automorphism of order three with nine
fix-points are rarer. If the $K3$-surface $X$ is algebraic, then its
Picard number is $\geq 19$. So in this case 
the surface $X$ and the covering surface $A$
can depend on at most one parameter.

Examples of abelian surfaces with an automorphism of
order three are given in [BH]: Each selfproduct $A=C \times C$,
with $C$ an elliptic curve, admits the automorphism
$$(x,y) \mapsto (-x,x-y).$$
It is shown in [BH] that the quotient $\bar{X}$
then is a double cover of the plane $\P_2$, branched over
the sextic $C^*$ dual to a plane cubic    
$C \subset \P_2$, a copy of the elliptic curve $C$. The
nine cusps of $\bar{X}$ of course ly over the nine cusps
of $C^*$.

By deformation theory of $K3$-surfaces, one may convince oneself, 
that there are also non-algebraic $K3$-surfaces with nine cusps.

\vv
Convention: Throughout this note the base field for algebraic 
varieties is $\C$.

\vV
\section{Cyclic triple covers of $K3$-surfaces}  

By a configuration of type $A_2$ on a smooth surface I
mean a pair $E,E'$ of smooth rational curves with
$E^2=(E')^2=-2, \quad E.E'=1$. Such a pair can be
contracted to a double point $A_2$ (a cusp).

{\bf Lemma 1.} {\em Let $X$ be a $K3$-surface carrying $p$ 
disjoint configurations
of type $A_2$, and $\bar{X}$ the surface obtained from $X$
by contracting them to cusps. If there is a 
smooth complex surface $Y$ and a triple
cover $Y \to \bar{X}$ branched (of order three) precisely over the
$p$ cusps, then
either $p=6$ with $Y$ a $K3$-surface or $p=9$ with
$Y$ a torus.}

Proof. First of all, $Y$ is k\"ahler: Indeed, $X$ is
k\"ahler by [S].
Blow up $X$ in the 
$p$ points, where the $p$ pairs of curves in the $A_2$ 
configurations meet. The resulting surface $\tilde{X}$
is k\"ahler by [B, Theoreme II 6]. Pull back the covering
$Y \to \bar{X}$ to a covering $\tilde{Y} \to \tilde{X}$.
Here $\tilde{Y}$ is k\"ahler, since it is a smooth 
surface in some $\P_1$-bundle over $\tilde{X}$, 
which is k\"ahler by [B, Theoreme principal II]. The surface
$Y$ is obtained from $\tilde{Y}$ by blowing down $(-1)$-curves,
so it is k\"ahler too by [F].

The canonical bundle of $Y$ admits a section
with  zeros at most in $p$ points. So there are no such zeros and
$K_Y$ is trivial. 

By the classification of surfaces
[BPV, p.188] the covering surface $Y$ therefore either is $K3$ with
$e(Y)=24$ or a torus with $e(Y)=0$.
The Euler number of $Y$ is computed in terms
of $p$ as
$$e(Y)= 3 \cdot e(\bar{X})- 2 \cdot p = 
3 \cdot(24-2p)-2 \cdot p = 72-8 \cdot p.$$
The possibilities are $p=6$ with $Y$ a $K3$-surface
and $p=9$ with $Y$ a complex torus.   \qed

\vv
Consider the $2p$ rational curves $E_i,E_i' \subset X$ 
forming the $p$ configurations of type $A_2$. The cyclic cover,
lifted to $X$, is branched along all of these curves of
order three. So there must be a divisor
$$\sum_{i=1}^p a_i \cdot E_i+a_i' \cdot E_i', \qquad 
a_i,a_i'=1 \mbox{ or }2,$$
divisible by three in $Pic(X)$. This implies that all
intersection numbers 
\begin{eqnarray*}
(a_i E_i+a_i'  E_i').E_i &=& -2 \cdot a_i+a_i' \\
(a_i E_i+a_i'  E_i').E_i' &=& a_i -2 \cdot a_i'
\end{eqnarray*} 
are divisible by three. Hence 
$$a_i=1 \Leftrightarrow a_i'=2 \quad
\mbox{ and } \quad a_i=2 \Leftrightarrow a_i'=1.$$ 

\vV
\section{The lattice generated by the nine cusps}

Here let $X$ be a $K3$-surface and $L=H^2(X,\Z)\simeq \Z^{22}$ its lattice
provided with the (unimodular) intersection form. Assume
that on $X$ there are nine disjoint $A_2$-configurations
$E_i,E_i', \, i=1,...,9$. Following [N] we denote by
$I \subset L$ the sublattice spanned (over $\Z$) by the
18 classes $[E_i], [E_i']$. Let me denote by 
$\bar{I} \subset L$ the {\em primitive} sublattice spanned
by these classes over $\Q$ and let me put
$$Q:=\bar{I}/I.$$
To study $Q$ we split $I$ in two sublattices by the
base change
$$ E_i,E_i' \quad  \mbox{ replaced by } \quad  E_i,F_i:=2 E_i+E_i'
\mbox{ for }i=1,...,9.$$
The essential point is that the intersection numbers
$$E_i.F_j=-3 \cdot \delta_{i,j}, \quad 
F_i.F_j= -6 \cdot \delta_{i,j}$$
are divisible by $3$.

{\bf Lemma 2.} {\em If a class
$$\sum_{i=1}^9 \epsilon_i [E_i] + \varphi_i[F_i],
\quad \epsilon_i, \varphi_i \in \Q,$$
belongs to $\bar{I}$ then
$$ \epsilon_i \in \Z, \quad 3 \cdot \varphi_i \in \Z.$$
In particular the order of the finite group $Q$ is
$|Q| =3^n$ for some $n \geq 0.$}

Proof. We just intersect the class with $E_k$ and $E_k'$ to find
$$ \begin{array}{c}\D
(\sum_{i=1}^9 \epsilon_i [E_i] + \varphi_i[F_i]).E_k = 
-2 \cdot \epsilon_k-3 \cdot \varphi_k \in \Z \\
  \\
\D
(\sum_{i=1}^9 \epsilon_i [E_i] + \varphi_i[F_i]).E_k' = 
\epsilon_k \in \Z. \\
\end{array}$$
This implies the assertion.   \qed

Lemma 2 shows in particular
$$\bar{I} = E + \bar{F}$$
with $E \subset I$ the lattice spanned by the classes $[E_i]$,
with $F \subset I$ the lattice spanned by the classes $[F_i]$,
and $\bar{F}$ the primitive sublattice of $L$ spanned over 
$\Q$ by $F$.

\vv
{\bf Lemma 3.} {\em The order $|Q|$ is $3^n$ with 
$n \geq 3$.}

Proof. Choose a system of $n$ generators for $Q$. They are
the residues of $n$ classes $q_1,...,q_n \in \bar{F}$.
The set of these $n$ classes can be extended to a 
$\Z$-basis 
$$q_1,...,q_n, f_{n+1},...,f_9$$
of $\bar{F}$. So, if $n \leq 2$, the lattice $\bar{F}$ has
a $\Z$-basis $q_1,q_2,f_3,...,f_9$ with $f_3,...,f_9$ integral
linear combinations of the classes $[F_1],...,[F_9]$. We extend
this basis to a $\Z$-basis of $\bar{I}$ with the classes 
$e_1=[E_1],...,e_9=[E_9]$, and to a basis of $L$ with some
classes $t_{19},...,t_{22}$. In the basis 
$$f_3,...,f_9,e_1,...,e_9,q_1,q_2,t_{19},...,t_{22}$$
the intersection matrix is
$$
\begin{array}{c|c|c}
\multicolumn{1}{c}{7} & \multicolumn{1}{c}{9} & \multicolumn{1}{c}{6} \\
(f_i.f_j) & (f_i.e_j) & * \\  \hline
(f_i.e_j) & (e_i.e_j) & * \\  \hline
  *     &   *     & * \\
\end{array}$$
Each summand in the Leibniz expansion of the determinant 
contains at least ten factors 
$$f_i.f_j, \quad  f_i.e_j \quad \mbox{ or } \quad e_i.e_j.$$
At most nine of them can be $e_i.e_j$.  
At least one of them must be a factor 
$ f_i.f_j$ or $f_i.e_j$ divisible by 3.
This shows that the determinant 
of the $22 \times 22$ intersection matrix is
divisible by $3$, a contradiction  with unimodularity. \qed

\vV
\section{The code of the nine cusps}

Each class in $\bar{I}$ is of the form
$$\sum_{i=1}^9 \epsilon_i \cdot E_i +\varphi_i \cdot F_i, \quad
\epsilon_i \in \Z, \, \varphi_i \in \frac{1}{3} \Z.$$
By sending
$$\varphi_i \mapsto \varphi_i \mbox{ mod } \Z$$
we identify $Q$ with an $\F_3$ sub-vector space of 
$\F_3^9$. By lemma 3 the sub-vector space $Q \subset \F_3^9$ has 
dimension $\geq 3$. In this section we want to identify this
sub-vector space.

In analogy with coding theory, we call each vector
$q = (q_i)_{i=1,..,9} \in Q$ a {\em word}, and the number
of its non-zero coefficients its {length} $|q|$.  By
lemma 1 all vectors $q \in Q$ have length
$|q|=0,6$ or $9$. As $dim_{\F_3}(Q) \geq 3$, 
the space $Q$ contains at least $3^3-3 =24$ words
of length $6$. 

\vv
We say that two words $q,q'$ overlap in $r$ places, if
there are precisely $r$ ciphers $i$ such that both 
coefficients $q_i$ and $q_i'$ are nonzero.
It is clear that any two nonzero words of length
six overlap in at least three places. If they overlap in six places,
they are linearly dependent: In fact, if $q+q' \not=0$, 
we have $q_i=q_i'$ for at least one $i$. Then $q+2q'$
has length $\leq 5$, hence $q+2q'=0$.

{\em Claim 1. Any two linearly independent vectors 
$q,q'$ of length six overlap in three or in four places. }

Proof. We have to exclude, that $q$ and $q'$
overlap in five places. Assume to the contrary that they do.
By rescaling the basis vectors of $\F^3$ we
may assume
$$q=(1,1,1,1,1,1,0,0,0)$$
and
$$q'=(0,q_2',q_3',q_4',q_5',q_6',q_7',0,0), \quad q_i'=1 \mbox{ or }2.$$
Since $q+q'$ again is a word of length six, w.l.o.g.
$$q'=(0,2,1,1,1,1,q_7',0,0).$$
Then 
$$q+2q' = (1,2,0,0,0,0,2q_7',0,0) \notin Q,$$
contradiction. \qed

\vv
Now, let me call the nine ciphers $1,...,9$ 'points'
and those triplets $\{i,j,k\}$ of ciphers 'lines', for which 
there is a word $q$ of length six with $q_i=q_j=q_k=0$.
As there are at least 24 words of length six, there
are at least twelve lines. As two linearly independent
words of length six overlap in four or three places,
two different lines intersect in one point, or not
at all (parallel lines). This allows to count the
number of lines:

{\em Claim 2. There are precisely $12$ lines, and therefore
the dimension of $Q$ is $n=3$.}

Proof. Through each point, there are at most four
distinct lines. So there are at most $9 \times 4 = 36$
incidences of lines with points. As on each line there
are three points, we have indeed at most $36/3=12$ lines. \qed

This proof shows in particular, that through each point
there are exactly four lines, or in other words: Each pair
of points lies on a (uniquely determined) line.

\vv
{\em Claim 3. For each line there are precisely two parallel lines. 
These two parallel lines do not intersect.}

Proof. Each line $L$ meets $3 \times 3=9$ other lines, hence there
are two lines $L',L''$ parallel to it. If $L'$ and $L''$ would meet
in a point, then through this point we would have five lines:
the two lines $L'$ and $L''$ and the three lines joining this point
with the three points on $L$. \qed

\vv
{\em Claim 4. The code $Q$ contains a word of length nine.}

Proof. Take three parallel lines $L,L',L''$  and two
words $q,q'$ vanishing on the lines $L,L'$ respectively. 
These two words $q$ and $q'$ overlap in precisely
three places (the points of $L''$). After
replacing $q$ by $2 \cdot q$ if necessary, we may assume 
$q_i=q_i'$ for one $i \in L''$. Then $q-q'$ is a word of
length six, i.e., $q_i=q_i'$ for all $i \in L''$.
So $q+q'$ is a word of length nine. \qed

\vv
Claim 4 proves the theorem from the introduction:
The existence of a word of length nine shows that
there is a linear combination
$$D:=\sum_{i=1}^9 \varphi_i F_i \in Pic(X) \quad
\mbox{with} \quad 0<\varphi_i<1, 3 \varphi_i \in \Z.$$
The divisor
$$3 \cdot D - \sum_{\varphi_i=2} 3 \cdot E_i'$$
contains all curves $E_i$ and $E_i', \, i=1,...9,$ with
multiplicity $1$ or $2$, and it is divisible by 3.

\vV
\section{The double cover branched over the dual cubic}

A smooth cubic $C \subset \P_2$ has nine flexes.
On the dual cubic $C^* \subset \P_2^*$ they yield nine cusps.
So the double cover $\bar{X} \to \P_2^*$ is an example
of a $K3$-surface with nine cusps. 
Here I want to understand the $3$-torsion property on $\bar{X}$
in terms of plane projective geometry, independently of the
theory in the preceding sections and of [BH].

The nine flexes of $C$ in a natural way have the structure
of an affine plane over $\F_3$. In fact, if $C$ is given
in Hesse normal form
$$x_0^3+x_1^3+x_2^3+3 \lambda x_0x_1x_2=0,$$
a transitive action of the vector space $\F_3^2$ on the curve $C$ and thus
on the set of its flexes is induced by the symmetries
$$\sigma:x_i \mapsto x_{i+1}, \quad 
\tau:x_i \mapsto \omega^i \cdot x_i, \quad
\omega \mbox{ a primitive third root of unity}.$$
Of course the 'lines' used in the preceding section
must be the lines in this affine plane. This
section will give a proof. 

Let me in this section denote by a line in the 
set of flexes, a line in the sense of the affine
structure just mentioned. 

The flexes are cut out on $C$ by the coordinate triangle
$x_0x_1x_2=0$. Two parallel lines are formed e.g. by the 
triplet of flexes $(0:1:-\omega^k)$ and the triplet 
$(1:0:-\omega^k)$. 
(All pairs of parallel lines are equivalent to this one, so
let us restrict our attention to this pair.)
The inflectional tangents there
are 
$$ -\lambda \omega^k \cdot x_0+ x_1+\omega^{2k} \cdot x_2=0
\quad \mbox{ and } \quad
x_0-\lambda \omega^k \cdot x_1+\omega^{2k} \cdot x_2=0.$$
The essential remark is, that they touch a nondegenerate 
conic, which in dual coordinates $(\xi_0:\xi_1:\xi_2)$
has the equation
$$\xi_0 \cdot \xi_1 + \lambda \xi_2^2=0.$$
(Of course, here we have to exclude $\lambda=0$, the case of the Fermat
cubic, where these triplets of inflectional tangents
are concurrent.) This implies that the corresponding six
cusps on the dual cubic $C^*$ in $\P_2^*$ are cut out by
the nondegenerate conic, whose equation was just
given. This conic intersects $C^*$ in each cusp
with multiplicity 2, so does not touch the
tangent of the cusp. 

Clearly, the inverse image of this conic on $\bar{X}$
decomposes into two smooth rational curves 
$\bar{R},\bar{R}' \subset \bar{X}$ passing through
our six distinguished cusps. Denote by $R,R'$ the proper
transforms of these curves on the smooth surface $X$.
A computation in local coordinates shows, that each
curve $R$ or $R'$ meets just one of the two rational curves
$E_i,E_i'$ from the $A_2$-configuration 
over each of the six distinguished cusps. Let
me call $E_i$ those curves which meet $R$, and $E_i'$
the curves intersecting $R', \, i=1,...,6$.

\vv
{\bf Lemma 4.} {\em For general choice of $\lambda$, the 
$K3$ surface $X$ has Picard number 19.}

By [PS, \S8] there is only a countable set of $K3$-surfaces with
Picard number 20. So, all we have to show is that the structure of $X$ indeed
varies with the elliptic curve $C$. In fact, a copy of $C$ (the proper
transform of the branch locus) lies on $X$, where it passes through
the intersection points in $E_i \cap E_i',\, i=1,...,6$. So, if
the structure of $X$ would not vary with $C$, we would have on
$X$ more than countably many elliptic pencils, 
a contradiction. \qed

This implies that $NS(X)$ is generated (over $\Q$) by the classes
of $E_i, E_i',\, i=1,...,9,$ and by the pullback $[H]$ of
the class of a line on $\P_2^*$. 
Now put
$$R-R' \sim \sum_{i=1}^9 (n_i \cdot [E_i]+n_i' \cdot [E_i'])
+n \cdot H, \; n_i,n_i',n \in \Q.$$
From
$$(R-R').H = (R-R').E_i = (R-R').E_i' =0
\mbox{ for } i=7,8,9$$
we conclude $n=n_i=n_i'=0$ for $i=7,8,9$. The other intersection numbers
are ($i=1,...,6$)
$$\begin{array}{rcccr}
1 &=& (R-R').E_i &=& -2n_i+n_i' \\
-1 &=& (R-R').E_i' &=& n_i-2n_i' \\
\end{array}$$
This implies
$$n_i' = -n_i = \frac{1}{3}.$$
We have shown that the class
$$\frac{1}{3} \sum_{i=1}^6 (E_i'-E_i) = R-R'$$
is integral. This is equivalent to the fact that the classes
$$\sum_{i=1}^6 (2 \cdot E_i+E_i')
\quad \mbox{ and } \quad \sum_{i=1}^6 (E_i+2 \cdot E_i')$$
are divisible by 3. 

Finally, we remark: For a 3-divisible set of six cusps on $\bar{X}$
the pattern, in which the curves $E_i$ and $E_i'$ organize themselves
into unprimed and primed ones, is given by their intersections
with $R$ or $R'$. 

\vV
\section{References}

\noindent
[BPV] Barth, W., Peters, C., Van de Van, A.: Compact complex
surfaces, Ergebnisse der Math. (3), 4, Springer (1984)

\noindent
[BH] Birkenhake C., Lange, H.: A family of abelian surfaces and
curves of genus four. manuscr. math. 85, 393-407 (1994)  

\noindent
[B] Blanchard, A.: Sur les varietes analytiques complexes.
Ann. Sci. ENS 73, 157-202 (1954)

\noindent
[F] Fujiki, A.: K\"{a}hlerian normal complex spaces. Tohoku Math. J.,
$2^{nd}$ series, 35, 101-118 (1983)

\noindent
[N] Nikulin, V.V.: On Kummer surfaces. Math. USSR Izv. 9, No 2, 261-275
(1975)

\noindent
[PS] Pjateckii-\v{S}apiro, I.I, \v{S}afarevi\v{c}, I.R.: A Torelli
theorem for algebraic surfaces of type $K3$. Izv. Akad. Nauk
SSSR, 35, 530-572 (1971) 

\noindent
[S] Siu, Y.T.: Every $K3$-surface is k\"ahler. Invent. math. 73,
139-150 (1983)

\end{document}